\documentclass[9pt,a4paper]{article}
\usepackage[utf8]{inputenc}
\usepackage[T1]{fontenc}
\usepackage{lmodern} 
\usepackage[francais,english]{babel}
\usepackage{amsmath,amsthm}
\usepackage{mathrsfs}
\usepackage{bbm} 
\usepackage{amssymb}
\usepackage{amsfonts}
\usepackage{pifont}
\usepackage{stmaryrd}
\usepackage{dsfont}
\usepackage{graphicx,pstricks}
\usepackage{boites,boites_exemples}
\usepackage{enumerate}
\usepackage{lscape}
\usepackage[all]{xy}
\usepackage{hyperref}
\input xy
\xyoption{all}
\newcommand{\be}{\begin{enumerate}}
\newcommand{\ee}{\end{enumerate}}
\newcommand{\R}{\mathbb{R}}

\newcommand{\C}{\textnormal{C}}

\newcommand{\W}{\textnormal{W}}

\renewcommand{\H}{\textnormal{H}}
\renewcommand{\L}{\textnormal{L}}

\newcommand{\T}{\mathbb{T}}

\renewcommand{\div}{\textnormal{div}}

\newcommand{\dd}{\mathrm{d}}
\newcommand{\na}{\nabla}
\newcommand{\pa}{\partial}

\newtheorem{thm}{Theorem}[section]

\newtheorem{rem}{Remark}[section]

\newcommand{\eps}{\varepsilon}

\numberwithin{equation}{section}

\author{
Daniel Han-Kwan
  \thanks{CMLS, CNRS, \'Ecole polytechnique, Institut Polytechnique de Paris, 91128 Palaiseau Cedex, France. Email: \textsf{daniel.han-kwan@polytechnique.edu}}
  \and
Mikaela Iacobelli
  \thanks{ETH Z\"urich, R\"amistrasse 101, 8092, Zurich. Email: \textsf{mikaela.iacobelli@math.ethz.ch}}
}

\title{From Newton's second law to Euler's equations of perfect fluids}

\begin{document}
\maketitle

\begin{abstract}
Vlasov equations can be formally derived from $N$-body dynamics in the mean-field limit. In some suitable singular limits, they may themselves converge to fluid dynamics equations.
Motivated by this heuristic, we introduce natural scalings under which the incompressible Euler equations can be rigorously derived from $N$-body dynamics with repulsive Coulomb interaction. Our analysis is based on the modulated energy methods of Brenier \cite{Br00} and Serfaty \cite{S_D18}.
\end{abstract}

\section{Introduction and main results}
Consider a system of $N$ indistinguishable particles evolving under the influence of binary interactions,
according to Newton's second law. For simplicity, we study periodic boundary conditions, which means that the phase space is $\T^d \times \R^d$, $d\geq 2$. We assume the interactions are described by an interaction force $-\na \Phi$ (with $\Phi$ which will be taken as the repulsive Coulomb kernel in the following). The phase space positions $(x_i(t),v_i(t))_{i=1}^N$ of the particles are then governed by the following system of ODEs:
\begin{equation}
\label{eq:ODE N 1}
\left\{
\begin{array}{l}
\displaystyle\frac{\dd}{\dd t} x_i = v_i,\medskip \\
\displaystyle \frac{\dd}{\dd t} v_i = -\lambda \sum_{\substack{j=1\\i\neq j}}^N \na \Phi (x_i -x_j).
\end{array}
\right.
\end{equation}
In these equations, the scaling parameter $\lambda$ is a function of $N$, which accounts for various intrinsic parameters of the system of particles.
In most physically relevant
cases, the number $N$ of particles is very large, as one can see by considering as an example the case of a rarefied gas or a plasma, for which $N$ is (at least) of order $10^{23}$. This means that in order to understand the dynamics of the system, even from the point of view of numerical analysis, it is unreasonable to determine the evolution of each individual particle exactly. 
Instead, it may be useful to replace such a microscopic model with a coarser description of the system.

A possible heuristic approach consists in considering the whole system as a fluid and to apply Newton's second law to each infinitesimal
volume element. In this way, one can formally recover equations from fluid mechanics,
such as the Euler or Navier-Stokes equations, which describe the evolution of macroscopic observable quantities associated with the
system, namely density, velocity, temperature, etc.

Kinetic equations, instead, describe the evolution of a $N$-body system at a mesoscopic
level, which is an intermediate viewpoint between the Newtonian dynamics and a macroscopic hydrodynamic model. From this point of view, the
evolution of the $N$-body system is described in a statistical manner by a single scalar quantity, namely the distribution function, which is governed by an (integro-)partial differential equation set on the single-particle phase space.

These three different levels of description of many-body systems can be translated into other mathematical
models. It is a fundamental problem of mathematical physics to understand how macroscopic behavior emerges from the underlying
dynamics of its constitutive particles.

The mean-field limit concerns the limit as the number of particles $N$ goes to infinity, from~\eqref{eq:ODE N 1} to a kinetic equation; in the case of Coulomb interaction, this corresponds to the derivation of the Vlasov-Poisson equation, a problem of tremendous importance on which a lot of progress has been lately made (for references on this topic, see in particular \cite{ Spohn91,HJ, G16, Hauray, Hauray-Jabin, J14, Laz, Laz-Pickl, S_D18}), although it is fair to say that a complete answer remains open in full generality.

This limiting procedure is usually performed in the so-called mean-field scaling, that is when the scaling parameter is $\lambda\sim 1/N.$ As explained in the introduction of the review paper \cite{J14} by Jabin, this choice is mathematically convenient, as it formally yields a force term of order $1$, assuming that each term of the sum is of order $1$.
However, it may be relevant to consider the case of scaling parameters $\lambda$ of order significantly larger than $1/N$, say $\lambda \sim 1/N^{\theta}$ with $\theta<1$, 
in which case the force term can in principle diverge, and one could expect a more singular behavior. We refer to this regime as a \emph{supercritical mean-field limit}; in full generality, the analysis is likely to be complex and heavily dependent on the structure of the interaction.
This question motivates the study of different singular limits than the usual mean-field limit in order to connect the different scales of description of the $N$-body system. 

To conclude, for what concerns the passage from the kinetic to the macroscopic description, one can rigorously derive the incompressible Euler equations from the Vlasov-Poisson equation through singular limits, namely, the \emph{quasineutral} \cite{Br00} and the \emph{gyrokinetic} \cite{GSR99} limits (which therefore may be thought to as hydrodynamic limits).

In a recent breakthrough \cite{S_D18}, Serfaty introduced a modulated energy method which allowed her to prove quantitative mean-field limits for a large class of first-order systems as well as for \eqref{eq:ODE N 1} in the \emph{monokinetic regime}, that is to say for velocities concentrating to the same vector field, leading to a derivation of the pressureless Euler-Poisson system (the later in an appendix in collaboration with Duerinckx). This method was subsequently applied in~\cite{GP} in the context of $N$-body quantum dynamics.
See also \cite{BJW,R2020} for further developments related to this approach.

Motivated by \cite{S_D18} and by the aforementioned problem of supercritical mean-field limit, we propose a natural scaling of the $N$-body problem (and variants thereof) which allows, in the monokinetic regime, to derive the \emph{incompressible Euler equations} (without passing via an intermediate kinetic model). Let us recall that they read as
\begin{equation}
\label{eq:Euler nD}
\pa_t u + u\cdot \na u = -\na p, \quad \div \,u =0.
\end{equation}
These equations describe the evolution of the velocity field $u=u(t,x)$ of a homogeneous incompressible inviscid (i.e. perfect) fluid in the absence of any other force than those arising from the incompressibility constraint. Here the expression $u\cdot \na u$ is the vector in $\R^d$ whose $i$-th component is given by $\sum_j u_j \partial_{j} u_i$. Hence, \eqref{eq:Euler nD} is a system of $(d+1)$ equations where the scalar function $p(t,x)$ is the pressure field associated to the incompressibility constraint.

We shall focus in this work on two regimes for systems of particles interacting through repulsive Coulomb potential, namely:

\begin{enumerate}

\item[(i)] In arbitrary dimension $d\geq 2$, we effectively investigate the {\bf supercritical mean-field limit}, which we may actually interpret as a combined mean-field and \emph{quasineutral} limit.
This question was recently studied in \cite{GPI18, GPI20}; in that papers, no monokinetic regime was required, and this allowed the authors to recover the so-called \emph{kinetic Euler equation} (loosely speaking, a kinetic generalization of \eqref{eq:Euler nD}, see e.g. \cite{BrCPAM}). However, at such level of generality, the scaling in \cite{GPI18, GPI20} is only slightly supercritical, in the sense that $\lambda \sim 1/ (N \log \log N)$.

\item[(ii)] In dimension $d=2$, we tackle the description of {\bf strongly magnetized plasmas}, leading to a combined mean-field and \emph{gyrokinetic} limit, for which we assume the presence of a fixed, constant \emph{strong magnetic field}, and we study the dynamics on a \emph{long time} scaling.
\end{enumerate}

\subsection{Supercritical mean-field limit, or the combined mean-field and quasineutral limits}
In this section we aim at studying~\eqref{eq:ODE N 1} for $\lambda \sim 1/N^\theta$, with $\theta <1$. Equivalently, we write $\lambda \sim \frac{1}{\eps^2}\frac{1}{N}$ with $\eps$ to be seen as a function of $N$. Focusing on the case of repulsive Coulomb interaction, this yields the $N$-body system:
\begin{equation}
\label{eq:ODE N}
\left\{
\begin{array}{l}
\displaystyle\frac{\dd}{\dd t} x_i = v_i,\medskip \\
\displaystyle \frac{\dd}{\dd t} v_i = -\frac{1}{\eps^2 N} \sum_{\substack{j=1\\i\neq j}}^N \na g (x_i -x_j).
\end{array}
\right.
\end{equation}
In this formulation, the parameter $\eps$ can be understood as the Debye length (that is the typical length scale of electrostatic interactions), and the limit as $\eps$ goes to zero is usually referred to as the quasineutral limit.
In this equation, $g$ stands for the Green function associated to the negative Laplacian on the torus, normalized to have zero mean:
\begin{equation}
\label{eq:g mean}
- \Delta g = \delta_0 - 1,\qquad\int_{\T^d} g \,\dd x=0,
\end{equation}
where $\delta_0$ stands for the Dirac mass with support $\{0\}$.
Thus $g \in C^\infty(\T^d \setminus \{0\}).$ Moreover, if we denote by $B_{1/4}(0)$
the ball of radius $1/4$ and centered at zero, then inside $B_{1/4}(0)$ the function $g$ can be decomposed as
\begin{equation} \label{def:G0}
g(x) = \begin{cases}
- \frac{1}{2 \pi} \log{|x|} + g_0(x) & d=2 \\
 \frac{1}{ |\mathbb{S}_{d-1}| |x|^{d-2}} + g_0(x) & d\ge3 ,
\end{cases}
\end{equation}
where $g_0 \in C^\infty(\overline{B_{1/4}(0)})$ is a smooth function, and $|\mathbb{S}_{d-1}|$ denotes the surface area of the unit sphere in dimension $d$ (see e.g. \cite[Lemma 3.4]{GPI-WP}).
We identify the $d$-dimensional torus $\T^d$ with the cube $\mathcal{Q}_d : = \left [- \frac{1}{2}, \frac{1}{2} \right ]^{d}$, with appropriate identifications of the boundary. 
In particular, since $x_i$ and $x_j$ denote points on the torus, by $x_i-x_j$ we mean the representative in $\mathcal{Q}_d$.
Note that the force term $\frac{1}{\eps^2 N} \sum_{\substack{j=1\\i\neq j}}^N \na g (x_i -x_j)$ can be equivalently written as
$$
\frac{1}{\eps^2}\int_{\T^d\setminus \{x_i\}} \nabla g(x_i-y)\, \dd\biggl( \frac{1}{N}\sum_{\substack{j=1\\i\neq j}}^N \delta_{x_j} \biggr)(y).
$$
We observe, since we deal with the repulsive Coulombian case, that the conservation of the total energy in the ODEs system implies both that velocities remain bounded and the particles cannot get too close to each other, uniformly in time.
In particular, as soon as the initial spatial positions of the $N$ particles are separated, \eqref{eq:ODE N} is a smooth ODE system which has a solution for all times.

The formal limit of \eqref{eq:ODE N} as $N\to \infty$ leads to the so-called Vlasov Poisson system.
The quasineutral limit (i.e., the limit $\eps\to 0$) of Vlasov-Poisson is a subject of its own which involves subtle stability questions; there have been several contributions about the rigorous passage to the limit under different hypotheses on the initial datum (see e.g., \cite{BG94,Grenier95,Grenier96,Br00,Mas01,HKH15,HKI2,HKI1,HKR}).
In particular, in \cite{Br00}, Brenier proves that in the monokinetic regime, the quasineutral limit of the Vlasov-Poisson equation leads to the incompressible Euler equations.
Therefore, taking the limit in \eqref{eq:ODE N}, first $N\to \infty$ then $\eps\to 0$ formally also leads, in the monokinetic regime, to the incompressible Euler equations. 
The purpose of this section is precisely to rigorously study the simultaneous limits $N\to \infty$ and $\eps\to 0$, seeing $\eps$ as a function of $N$.

Inspired by \cite{Br00} and the appendix of \cite{S_D18}, we introduce the following modulated energy :
\begin{equation}
\label{eq:H1}
\begin{aligned}
H_{N,\eps} (t) &= \frac{1}{2 N} \sum_{i=1}^N |u(t,x_i) -v_i|^2 \\
&+ \frac{1}{2\eps^2} \iint_{\T^d \times \T^d \setminus \triangle} g(x-y)\, \dd\biggl( \frac{1}{N}\sum_{i=1}^N \delta_{x_i} - 1- \eps^2 \mathfrak{U} \biggr)(x)\, \dd\biggl( \frac{1}{N}\sum_{i=1}^N \delta_{x_i} - 1- \eps^2 \mathfrak{U} \biggr)(y)\\
&=: H_{N,\eps}^{(1)}(t) + H_{N,\eps}^{(2)} (t),
\end{aligned}
\end{equation}
where $\triangle:=\{(x,y) \in \mathbb T^d\times \mathbb T^d\,:\,x=y\}$ is the diagonal in $\mathbb T^d\times \mathbb T^d$, $u$ satisfies the incompressible Euler equations \eqref{eq:Euler nD}
and we set 
\begin{equation}
\label{eq:U}
\mathfrak{U} = \div \, \div (u \otimes u) = \sum_{i,j} \pa_{x_i} u_j\, \pa_{x_j} u_i.
\end{equation}
Taking the divergence in \eqref{eq:Euler nD}, we have
\begin{equation}
\label{eq:p U}
-\Delta p=\mathfrak{U}.
\end{equation}
It turns out that $H_{N,\eps}$ can be seen as the $N$-particle version of the modulated energy introduced by Brenier \cite{Br00} to study the quasineutral limit of the Vlasov-Poisson system in the monokinetic regime. Note moreover that we need to include an additional corrector of order $\eps^2$ in $H_{N,\eps}^{(2)}$ that will play an important role in the analysis.

Our main result reads as follows:
\begin{thm}\label{thm:1}
Let $u \in L^\infty([0,T],C^{1,\alpha}(\mathbb T^d))$ be a classical solution of the incompressible Euler equations \eqref{eq:Euler nD}.
Let $(x_i(0),v_i(0))_{i=1}^N$ be a family of points in $\mathbb T^d\times \mathbb R^d$ such that $x_i(0) \neq x_j(0)$ for $i\neq j$, 
and let $(x_i(t),v_i(t))_{i=1}^N$ solve the ODE system \eqref{eq:ODE N}. 
Suppose that
$\eps=\eps_N \to 0$ in a sufficiently slow way so that $H_{N,\eps}(0)\to 0$ and $\eps N^{\frac{1}{d(d+1)}} \to \infty$ as $N\to \infty$. 
Then
$$
\sup_{t \in [0,T]}H_{N,\eps}(t) \to 0\qquad \text{as }N\to \infty.
$$
In particular, for any $t \in [0,T],$ the empirical measure $\frac{1}{N}\sum_{i=1}^N\delta_{(x,v)=(x_i(t),v_i(t))}$ associated to the ODE system \eqref{eq:ODE N}
converges weakly-$\star$ in $\mathbb T^d\times \mathbb R^d$ to the measure $\dd x \otimes \delta_{v=u(t,x)}$.
\end{thm}

\begin{rem}
\label{rmk:thm1}
Some general comments are at order.
\begin{itemize}
\item[(i)]
The existence of a solution $u \in L^\infty([0,T],C^{1,\alpha}(\mathbb T^d))$ of \eqref{eq:Euler nD}
is always guaranteed (at least on some finite time interval depending on the initial datum, and in dimension $d=2$ one can even take $T=+\infty$) provided that $u_0 \in C^{1,\alpha}(\mathbb T^d)$,
see for instance \cite{BM}.
\item[(ii)]
The fact that the convergence $H_{N,\eps}(t) \to 0$ implies the weak-$\star$ convergence of the empirical measure has been proved in \cite{S_D18}.

\item[(iii)] The initial convergence $H_{N,\eps}(0)\to 0$ is a mild assumption, see Section~\ref{sec-init}.

\item[(iv)] The convergence of $H_{N,\eps}(t)$ is actually quantitative, see~\eqref{quanti} below.

\item[(v)] This result corresponds to the supercritical mean-field limit with $\lambda = \frac{1}{\eps^2 N} \sim \frac{1}{N^\theta}$ and $\theta \in \biggl(1- \frac{2}{d(d+1)},1\biggr)$.

 \item[(vi)] We do not know if the constraint $\eps N^{\frac{1}{d(d+1)}} \to \infty$ can be improved. We may note however, that it shows up in two seemingly unrelated arguments: once for the statistical relevance of the assumption that the initial modulated energy converges to $0$, and once when applying the stability results of \cite{S_D18}. Certainly, when $\eps$ is too small, the limit behavior should be utterly different from the incompressible Euler dynamics that we have obtained here.

\end{itemize}

\end{rem}

\begin{rem}[A Rescaling] Setting $w_j := \eps {v_j}$, \eqref{eq:ODE N} rewrites as
\begin{equation*}
\left\{
\begin{array}{l}
\eps \displaystyle \frac{\dd}{\dd t} x_i = w_i, \\
\eps \displaystyle \frac{\dd}{\dd t} w_i = - \frac{1}{N} \sum_{\substack{i,j=1\\i\neq j}}^N \na g (x_i -x_j).
\end{array}
\right.
\end{equation*}
for which the limit $N\to +\infty$ can therefore be interpreted as a \emph{long time} mean-field limit. Reading Theorem~\ref{thm:1} in the variables $(x_i, w_i)$, one obtains the long time stability (namely for times of order $1/\eps\sim N^\alpha $ for $\alpha < \frac{1}{d(d+1)}$) of 
the Dirac mass $\delta_{v=0}$. This could be seen in analogy to the results of Caglioti-Rousset \cite{CR07,CR08} who proved similar stability results for Penrose stable equilibria for smooth interaction kernels
(see also \cite{HKN} for related instability results).

\end{rem}

\subsection{Combined mean-field and gyrokinetic limit}

Let $\Omega$ be equal either to $\T^2$ or to $\R^2$,
and consider the $N$-body dynamics on $\Omega\times \R^2$. We assume here particles are also subject to a Lorentz force because of the presence of a fixed, constant, strong magnetic field.
As is usual in the literature, the word gyrokinetic will precisely refer to this strong magnetic field regime.
 With an appropriate long time scaling, the system of ODEs then reads
\begin{equation}
\label{eq:ODE 2}
\left\{
\begin{array}{l}
\displaystyle\eps \frac{\dd}{\dd t} x_i = v_i,\medskip\\
\displaystyle\eps \frac{\dd}{\dd t} v_i = - \frac{1}{N} \sum_{\substack{i,j=1\\i\neq j}}^N \na g (x_i -x_j) + \frac{v_i^\perp}{\eps}.
\end{array}
\right.
\end{equation}
Here, given a vector $v=(v_1,v_2)$ we use the notation $v^\perp=(-v_2,v_1)$.
Again, $g$ stands for the Green kernel associated to the negative Laplacian on $\Omega$. The purpose of this section is to investigate the asymptotics of \eqref{eq:ODE 2} in the simultaneous limit as $\eps\to 0$, $N\to +\infty$, which results in a combined mean-field and gyrokinetic limit.
This problem corresponds to a $N$-particle system version of the gyrokinetic limit for the Vlasov-Poisson equation studied by Golse and Saint-Raymond \cite{GSR99,SR02} employing compactness methods (see also the recent works \cite{Miot16,Miot19}). The limit equations are the incompressible Euler equations. Consequently, these are also the expected limit of \eqref{eq:ODE 2}.

\begin{rem}[A Rescaling] Setting $w_j := \frac{v_j}{\eps}$, \eqref{eq:ODE 2} rewrites as
\begin{equation*}
\left\{
\begin{array}{l}
\displaystyle \frac{\dd}{\dd t} x_i = w_i, \\
\displaystyle \frac{\dd}{\dd t} w_i = - \frac{1}{\eps^2} \Biggl(\frac{1}{N} \sum_{\substack{i,j=1\\i\neq j}}^N \na g (x_i -x_j) + {w_i^\perp}\Biggr),
\end{array}
\right.
\end{equation*}
which is closer to~\eqref{eq:ODE N}.
This turns out to be the $N$-particle system version of the gyrokinetic limit studied by Brenier in \cite{Br00}. 
\end{rem}

Again inspired by \cite{Br00} and \cite{S_D18}, we introduce the following modulated energy:
\begin{equation}
\label{eq:H2}
\begin{aligned}
\mathscr{H}_{N,\eps} (t) &= \frac{1}{2 N}  \sum_{i=1}^N \biggl|u(t,x_i) -\frac{v_i}{\eps} \biggr|^2 \\
&+ \frac{1}{2\eps^2} \int_{\Omega \times \Omega \setminus \Delta} g(x-y)\, \dd\biggl( \frac{1}{N}\sum_{i=1}^N \delta_{x_i} -\omega - \eps^2 \mathfrak{U}\biggr)(x)\, \dd\biggl( \frac{1}{N}\sum_{i=1}^N \delta_{x_i} -\omega - \eps^2 \mathfrak{U}\biggr)(y) \\
&=: \mathscr{H}_{N,\eps}^{(1)}(t) + \mathscr{H}_{N,\eps}^{(2)} (t),
\end{aligned}
\end{equation}
where 
\begin{equation}
\label{eq:coupling}
u = \na^\perp \psi, \quad \omega = \Delta \psi,
\end{equation}
and $(\omega,u)$ solves
\begin{equation}
\label{eq:Euler 2D}
\pa_t \omega + \div(u\,\omega) = 0,
\end{equation}
which corresponds to 2D Euler in vorticity formulation. 
In particular $u$ solves \eqref{eq:Euler nD}.
Also, as before, we set
$\mathfrak{U} = \div \, \div (u \otimes u).$ 

Our result reads as follows:
\begin{thm}\label{thm:2}
Let $\Omega=\T^2$ or $\Omega=\R^2$.
Let $(\omega,u)  \in L^\infty([0,T],L^1(\Omega)\cap C^{0,\alpha}(\Omega))\times L^\infty([0,T],C^{1,\alpha}(\Omega)\cap \dot H^1(\Omega))$ be a classical solution of the 2D Euler equation in vorticity formulation \eqref{eq:coupling}-\eqref{eq:Euler 2D}.
Let $(x_i(0),v_i(0))_{i=1}^N$ be a family of points in $\Omega\times \mathbb R^2$ such that $x_i(0) \neq x_j(0)$ for $i\neq j$ and let $(x_i(t),v_i(t))_{i=1}^N$ solve the ODE system \eqref{eq:ODE 2}. 
Suppose that 
$\eps=\eps_N \to 0$ in a sufficiently slow way so that $\mathscr H_{N,\eps}(0)\to 0$ and $\eps N^{\frac{1}{6}} \to \infty$ as $N\to \infty.$
Then
$$
\sup_{t \in [0,T]}\mathscr H_{N,\eps}(t) \to 0\qquad \text{as }N\to \infty.
$$
In particular, the rescaled empirical measure $\frac{1}{N}\sum_{i=1}^N\delta_{(x,v)=(x_i(t),\eps^{-1}v_i(t))}$ associated to the ODE system \eqref{eq:ODE 2}
converges weakly-$\star$ in $\Omega\times \mathbb R^2$ to the measure $\omega \otimes \delta_{v=u(t,x)}$.
\end{thm}

\begin{rem}
A combined gyrokinetic and quasineutral limit has been considered in \cite{GSR03}; likely, such a regime can also be studied in the context of the mean-field limit. 
\end{rem}

The proofs of Theorems~\ref{thm:1} and~\ref{thm:2} are given in Sections~\ref{sec2} and~\ref{sec3}.\\
 
\subsection{On the statistical relevance of the initial convergence}
\label{sec-init}
In this section, we explain why considering initial configurations $(x_i(0),v_i(0))_{i=1}^N$ such that $H_{N,\eps} (0)\to 0$ (or $\mathscr{H}_{N,\eps} (0)\to 0$) as $\eps \to 0$ and $N \to \infty$ is statistically relevant in the regime $\eps N^{\frac{1}{d(d+1)}} \to \infty$.

Let $\mu \in L^1 \cap C^{0,\alpha}(\T^d)$. To simplify, we assume $d \geq 3$ (for $d=2$ some logarithmic corrections appear). We study the following functional 
\begin{align*}
 &\frac{1}{\eps^2}\mathbb{E} \biggl(\iint_{\Omega \times \Omega \setminus \Delta} g(x-y)\, \dd\biggl( \frac{1}{N}\sum_{i=1}^N \delta_{x_i} -\mu\biggr)(x)\, \dd\biggl( \frac{1}{N}\sum_{i=1}^N \delta_{x_i} - \mu \biggr)(y)\biggr) \\
 &= \frac{1}{\eps^2 N^2} \mathbb{E} \biggl( \sum_{i\neq j} g(x_i-x_j) \biggr)+ \frac{1}{\eps^2}\iint_{\Omega \times \Omega \setminus \Delta} g(x-y)\, \dd \mu (x)\, \dd \mu(y) \\
 &\quad - \frac{2}{\eps^2} \mathbb{E} \biggl(\int_{\Omega} (-\Delta)^{-1}\mu(x ) \dd\biggl( \frac{1}{N}\sum_{i=1}^N \delta_{x_i} \biggr)(x)\biggr),
\end{align*}
when the $(x_i)$ are $\mu$-distributed independent random variables.
It is quite straightforward to study the limit of the last term, as $ (-\Delta)^{-1}\mu$ is smooth:
$$
\frac{1}{\eps^2}\biggl|\mathbb{E} \biggl( \int_{\Omega} (-\Delta)^{-1}\mu(x ) \dd\biggl( \frac{1}{N}\sum_{i=1}^N \delta_{x_i} -\mu\biggr)(x) \biggr) \biggr| \lesssim \frac{1}{\eps^2}\mathbb{E} \biggl(W_1 \biggl( \frac{1}{N}\sum_{i=1}^N \delta_{x_i} , \mu \biggr)\biggr)\lesssim \frac{1}{\eps^2 N^{1/d}},
$$
according to \cite[Theorem 1]{FG}. Therefore to prove what we want, it is sufficient to ensure
$$
\frac{1}{\eps^2} \biggl|\mathbb{E} \biggl( \frac{1}{N^2} \sum_{i\neq j} g(x_i-x_j) - \iint_{\Omega \times \Omega \setminus \Delta} g(x-y)\, \dd \mu (x)\, \dd \mu(y) \biggr)  \biggr| \to 0.
$$
The idea is as follows. Let $\varphi \in C^\infty_c(B_1(0))$, with $0\leq \varphi\leq 1$ and $\varphi|_{B_{1/2}(0)}\equiv 1$. Given $\eta>0$ small, we write
$$
g_<(x) := g(x) \varphi\Bigl(\frac{x}{\eta}\Bigr), \quad g_> := g - g_<,
$$
and
$$
\frac{1}{N^2} \sum_{i\neq j} g(x_i-x_j) = \frac{1}{N^2} \sum_{i\neq j} g_<(x_i-x_j) + \frac{1}{N^2} \sum_{i\neq j} g_>(x_i-x_j) =: S_1 + S_2.
$$
The term $S_2$ is easy to study as $g_>$ is a smooth kernel:
$$
\begin{aligned}
\mathbb{E} \biggl| \frac{1}{N^2} \sum_{i\neq j} g_>(x_i-x_j) - \iint_{\Omega \times \Omega \setminus \Delta} g_>(x-y)\, \dd \mu (x)\, \dd \mu(y)  \biggr| &\lesssim \mathbb{E} \biggl(W_1 \biggl( \frac{1}{N}\sum_{i=1}^N \delta_{x_i} , \mu \biggr)\biggr) \| g_>\|_{W^{1,\infty}(\Omega)} \\
 &\leq \mathbb{E} \biggl(W_1 \biggl( \frac{1}{N}\sum_{i=1}^N \delta_{x_i} , \mu  \biggr) \biggr) \frac{1}{\eta^{d-1}} \\
 &\lesssim \frac{1}{N^{\frac{1}{d}}} \frac{1}{\eta^{d-1}},
 \end{aligned}
 $$
 On the other hand, we can write
 $$
 \mathbb{E} (S_1) = \frac{N-1}{N} \mathbb{E} \biggl(g_<(x_1 - x_2) \biggr)=  \frac{N-1}{N} \iint_{\Omega \times \Omega } g_<(x-y)\, \dd \mu (x)\, \dd \mu(y).
 $$
 We deduce by simple convolution estimates that
 $$
 | \mathbb{E} (S_1)| \lesssim \| g_<\|_{L^1} \| \mu \|_{L^\infty} \| \mu\|_{L^1} \lesssim \eta^2.
 $$
Optimizing in $\eta$ so that
$$
 \frac{1}{N^{\frac{1}{d}}} \frac{1}{\eta^{d-1}} = \eta^2,
$$
we conclude that
$$
\frac{1}{\eps^2}\biggl|\mathbb{E} \biggl( \frac{1}{N^2} \sum_{i\neq j} g(x_i-x_j) - \iint_{\Omega \times \Omega \setminus \Delta} g(x-y)\, \dd \mu (x)\, \dd \mu(y) \biggr)  \biggr| \lesssim \frac{1}{\eps^2 N^{\frac{2}{d(d+1)}}}.
$$
This proves that
$$
\frac{1}{\eps^2}\mathbb{E} \biggl(\iint_{\Omega \times \Omega \setminus \Delta} g(x-y)\, \dd\biggl( \frac{1}{N}\sum_{i=1}^N \delta_{x_i} -\mu\biggr)(x)\, \dd\biggl( \frac{1}{N}\sum_{i=1}^N \delta_{x_i} - \mu \biggr)(y)\biggr) \to 0
$$
as soon as $\eps N^{\frac{1}{d(d+1)}} \to \infty$. This is exactly the same constraint on $\eps$ as in the statement of Theorem~\ref{thm:1} (and similarly for Theorem~\ref{thm:2}, although in this case one needs to take care of the logarithmic corrections in dimension $2$).

We apply these considerations to $\mu \equiv 1$ (and the $(x_i(0))$ accordingly randomly picked) for what concerns Theorem~\ref{thm:1}.
Then one can choose the $(v_i(0))$ such that $|v_i(0)-u(0,x_i(0))|\leq \eta_N$ for all $i$, where $\eta_N \to 0$ as $N \to \infty$. The same discussion applies to Theorem~\ref{thm:2} as well with $\mu\equiv \omega$ (note that, in order to apply \cite[Theorem 1]{FG} in the case $\Omega =\R^2$, one needs to impose $\int_{\R^2 } \omega(0) |x|^q \, \dd x < +\infty$ for some $q>2$).
\\

 \noindent
 {\it Acknowledgements:} We wish to thank Francois Golse for inspiring discussions about this topic, and Pierre-Emmanuel Jabin for describing to us the typicality argument in Section~\ref{sec-init}. We also thank Megan Griffin-Pickering for her useful feedback on a preliminary version of this manuscript.
This work has been initiated during the INdAM workshop \emph{Recent advances in kinetic equations and applications} held at the University of Rome {Sapienza} in November $2019$.
DHK acknowledges the partial support of the grant ANR-19-CE40-0004.

\section{Proof of Theorem~\ref{thm:1}}
\label{sec2}
Recalling \eqref{eq:H1}, 
let us compute the derivative in time of the functional $H_{N,\eps}$.
Using \eqref{eq:ODE N}, by straightforward computations, we get the following formula:
\begin{equation*}
\begin{aligned}
\frac{\dd}{\dd t} H_{N,\eps}  (t)
&=\frac{1}{N} \sum_{i=1}^N  (u(t,x_i) - v_i) \cdot \biggl( \pa_t u(t,x_i) +  v_i \cdot \na u(t,x_i) + \frac{1}{\eps^2 N} \sum_{j,j\neq i} \na g (x_i -x_j)  \biggr) \\
&\quad+    \frac{1}{\eps^2N^2}\sum_{\substack{i,j=1\\i\neq j}}^N v_i \cdot \na g( x_i - x_j) +  \iint_{\T^d \times \T^d \setminus \triangle} g(x-y)   \,(1+\eps^2\mathfrak{U} (t,x))\, \pa_t \mathfrak{U} (t,y) \, \dd x \,\dd y  
\\
&\quad-   \frac{1}{N} \sum_{i=1}^N \int_{\T^d} v_i\cdot \na g (x_i - y) \, \mathfrak{U} (t,y)\, \dd y -   \frac{1}{N} \sum_{i=1}^N \int_{\T^d}  g (x_i - y)  \,\pa_t \mathfrak{U} (t,y) \,\dd y,
\end{aligned}
\end{equation*}
that is to say
\begin{equation}
\label{eq:dt H}
\begin{aligned}
\frac{\dd}{\dd t} H_{N,\eps}  (t)
&=\underbrace{\frac{1}{N} \sum_{i=1}^N  (u(t,x_i) - v_i) \cdot \biggl( \pa_t u(t,x_i) +  v_i \cdot \na u(t,x_i)\biggr)  -   \frac{1}{N} \sum_{i=1}^N \int_{\T^d} v_i\cdot \na g (x_i - y) \, \mathfrak{U} (t,y)\, \dd y}_{I} \\
&\quad+ \underbrace{\frac{1}{\eps^2 N} \sum_{\substack{i,j=1\\i\neq j}}^N u(t,x_i)\cdot \na g (x_i -x_j)}_{II}  \\
&\quad + \underbrace{\iint_{\T^d \times \T^d \setminus \triangle} g(x-y)   \,(1+\eps^2\mathfrak{U} (t,x))\, \pa_t \mathfrak{U} (t,y) \, \dd x \,\dd y  
-   \frac{1}{N} \sum_{i=1}^N \int_{\T^d}  g (x_i - y)  \,\pa_t \mathfrak{U} (t,y) \,\dd y}_{III}.
\end{aligned}
\end{equation}
We now need to properly manipulate each term in the right hand side.
We start with $I$.
First of all, 
since $u$ solves \eqref{eq:Euler nD}, we have
\begin{equation}
\label{eq:1}
\pa_t u(t,x_i) + v_i \cdot \na u(t,x_i) = (v_i- u(t,x_i)) \cdot \na u (x_i) - \na p (t,x_i).
\end{equation}
Also, using \eqref{eq:p U} we get
\begin{equation}
\label{eq:2}
 \sum_{i=1}^N  \int_{\T^d}  v_i\cdot \na g (x_i - y) \,  \mathfrak{U} (t,y)\, \dd y =  -\sum_{i=1}^N v_i \cdot \int_{\T^d}\na g (x_i - y) \,\Delta p(t,y)\, \dd y =   \sum_{i=1}^N v_i  \cdot  \na p (t,x_i ).
\end{equation}
Therefore, combining \eqref{eq:1} and \eqref{eq:2}, we get
\begin{equation}
\label{eq:12}
I=-\frac{1}{N} \sum_{i=1}^N (u(t,x_i) -v_i)\cdot \left[(u(t,x_i)-v_i) \cdot  \na u(t,x_i)\right]-\frac{1}{N} \sum_{i=1}^N  u(t,x_i)\cdot \nabla p(t,x_i).
\end{equation}
We now observe the following identity: 
$$
\begin{aligned}
&\frac{1}{2\eps^2} \iint_{\T^d \times \T^d \setminus \triangle} (u(t,x) - u(t,y)) \cdot \na g (x-y)\, \dd\biggl(  \frac{1}{N}\sum_{i=1}^N \delta_{x_i} - 1- \eps^2 \mathfrak{U} \biggr)(x)\, \dd\biggl(  \frac{1}{N}\sum_{i=1}^N \delta_{x_i} - 1- \eps^2 \mathfrak{U} \biggr)(y) \\
&\quad = \frac{1}{\eps^2N^2}   \sum_{\substack{i,j=1\\i\neq j}}^N u(t,x_i) \cdot \na g (x_i -x_j) - \frac{1}{N} \sum_{i=1}^N u(t,x_i)  \cdot  \na p (t,x_i )  \\
&\quad\quad +  \frac{1}{N} \sum_{i=1}^N \int_{\T^d} u(t,y)\cdot \na g(x_i-y)  \,\mathfrak{U} (t,y) \, \dd y - \iint_{\T^d \times \T^d \setminus \triangle } \na g (x-y) \cdot  u(t,y)\, \mathfrak{U} (t,y) \, \dd y \,\dd x. \\
&\quad \quad-  \eps^2  \int_{\T^d}   u(t,y)\cdot \na g(x-y)  \,\mathfrak{U} (t,x)\, \mathfrak{U} (t,y) \, \dd y \, \dd x\\
&\quad= \frac{1}{\eps^2N^2}   \sum_{\substack{i,j=1\\i\neq j}}^N u(t,x_i) \cdot \na g (x_i -x_j) - \frac{1}{N} \sum_{i=1}^N u(t,x_i)  \cdot  \na p (t,x_i )  \\
&\quad\quad+ \int_{\T^d\times \T^d \setminus \Delta} u(t,y)\cdot \na g(x-y)   \, \mathfrak{U} (t,y) \,   \dd y \,  \dd\biggl(  \frac{1}{N}\sum_{i=1}^N \delta_{x_i} - 1- \eps^2 \mathfrak{U} \biggr)(x) \\
&\quad =\frac{1}{\eps^2N^2}   \sum_{\substack{i,j=1\\i\neq j}}^N u(t,x_i) \cdot \na g (x_i -x_j) - \frac{1}{N} \sum_{i=1}^N u(t,x_i)  \cdot  \na p (t,x_i )  \\
&\quad\quad+ \int_{\T^d} \na (-\Delta)^{-1} \left( u\,  \mathfrak{U}\right)(x)\, \dd\biggl(  \frac{1}{N}\sum_{i=1}^N \delta_{x_i} - 1- \eps^2 \mathfrak{U} \biggr)(x).
\end{aligned}
$$
Thus, combining this  with \eqref{eq:12}, we obtain
\begin{equation}
\label{eq:13}
\begin{aligned}
I+II&= -\frac{1}{N} \sum_{i=1}^N (u(t,x_i) -v_i)\cdot \left[(u(t,x_i)-v_i) \cdot  \na u(t,x_i)\right]  \\
&+ \frac{1}{2\eps^2} \iint_{\T^d \times \T^d \setminus \triangle} (u(t,x) - u(t,y)) \cdot \na g (x-y)\, \dd\biggl(  \frac{1}{N}\sum_{i=1}^N \delta_{x_i} - 1- \eps^2 \mathfrak{U} \biggr)(x)\, \dd\biggl(  \frac{1}{N}\sum_{i=1}^N \delta_{x_i} - 1- \eps^2 \mathfrak{U} \biggr)(y) \\
&- \int_{\T^d} \na (-\Delta)^{-1} \left( u \,\mathfrak{U}\right)(x)\, \dd\biggl(  \frac{1}{N}\sum_{i=1}^N \delta_{x_i} - 1- \eps^2 \mathfrak{U} \biggr)(x).
\end{aligned}
\end{equation}
Finally, using again \eqref{eq:p U}, we get
\begin{equation}
\label{eq:5}
\begin{aligned}
III= \int_{\T^d}  \pa_t p (t,x)\, \dd\biggl(  \frac{1}{N}\sum_{i=1}^N \delta_{x_i} - 1- \eps^2 \mathfrak{U} \biggr)(x).
\end{aligned}
\end{equation}
Hence, combining \eqref{eq:dt H}, \eqref{eq:13}, and \eqref{eq:5}, we can rewrite
\begin{equation}
\begin{aligned}
\label{eq:dt H 2}
\frac{\dd}{\dd t} H_{N,\eps} &= -\frac{1}{N} \sum_{i=1}^N (u(t,x_i) -v_i)\cdot \left[(u(t,x_i)-v_i) \cdot  \na u(t,x_i)\right]  \\
&+ \frac{1}{2\eps^2} \iint_{\T^d \times \T^d \setminus \triangle} (u(t,x) - u(t,y)) \cdot \na g (x-y)\, \dd\biggl(  \frac{1}{N}\sum_{i=1}^N \delta_{x_i} - 1- \eps^2 \mathfrak{U} \biggr)(x)\, \dd\biggl(  \frac{1}{N}\sum_{i=1}^N \delta_{x_i} - 1- \eps^2 \mathfrak{U} \biggr)(y) \\
&- \int_{\T^d} \na (-\Delta)^{-1} \left( u \,\mathfrak{U}\right)(x)\, \dd\biggl(  \frac{1}{N}\sum_{i=1}^N \delta_{x_i} - 1- \eps^2 \mathfrak{U} \biggr)(x)  \\
&- \int_{\T^d}  \pa_t p (t,x)\, \dd\biggl(  \frac{1}{N}\sum_{i=1}^N \delta_{x_i} - 1- \eps^2 \mathfrak{U} \biggr)(x)  =: \sum_{i=1}^4 I_i.
\end{aligned}
\end{equation}
We now need to bound the four term $I_i$, $i=1,2,3,4.$

For the first integral, we have (recall \eqref{eq:H1})
\begin{equation}\label{eq:I1}
|I_1| \leq \| \na u\|_{\L^\infty(\T^d)} \frac{1}{N} \sum_{i=1}^N |u(t,x_i) -v_i|^2 \lesssim \| \na u\|_{\L^\infty(\T^d)}   H_{N,\eps}^{(1)}.
\end{equation}
For the other terms, we reply on two deep results from \cite{S_D18}.
More precisely, applying \cite[Proposition 1.1]{S_D18} 
with $s=d-2$ and normalizing both the left and right hand side appearing in the proposition by $\eps^{-2}N^{-2}$, we get
\begin{equation}\label{eq:I2}
\begin{aligned}
|I_2| \lesssim &\| \na u\|_{\L^\infty(\T^d)}  \biggl( H_{N,\eps}^{(2)} + \frac{1}{\eps^2} (1+ \eps^2 \|\mathfrak{U} \|_{\L^\infty(\T^d)}) \frac{1}{N^{\frac{2}{d(d+1)}}} + 1_{\{d=2\}} \frac{1}{\eps^2}  \frac{\log N}{N} \biggr) \\
&+ \frac{1}{\eps^2}\biggl( \| u \|_{\W^{1,\infty}(\T^d)} (1+ \eps^2 \| \mathfrak{U}\|_{\L^\infty(\T^d)} ) \frac{1}{N^{\frac{1}{d}}} + \| \na u \|_{\L^\infty(\T^d)} (1+ \eps^2  \| \mathfrak{U}\|_{\L^\infty(\T^d)} )  \frac{1}{N^{\frac{2}{d}}}\biggr).
\end{aligned}
\end{equation}
On the other hand, \cite[Proposition 3.6]{S_D18} rescaled by $N^{-1}$,
 we deduce that  
for all $\alpha \in (0,1]$ there exists $\lambda=\lambda(\alpha,d)>0$ such that 
\begin{equation}\label{eq:I34}
\begin{aligned}
|I_3| &+|I_4| \lesssim  \biggl(\|\na (-\Delta)^{-1} \left( u \,\mathfrak{U}\right)\|_{\C^{0,\alpha}(\T^d)} +\|\pa_t p \|_{\C^{0,\alpha}(\T^d)} \biggr) \frac{1}{N^{\frac{\lambda}{d}}} \\
&+\frac{1}{N}  \biggl(\|\na (-\Delta)^{-1} \left( u \,\mathfrak{U}\right)\|_{\dot\H^{1}(\T^d)} +\|\pa_t p \|_{\dot\H^{1}(\T^d)} \biggr) \Bigg(\eps^2 N^2 H_{N,\eps}^{(2)} +  (1+ \eps^2 \| \mathfrak{U}\|_{\L^\infty(\T^d)} ) N^{2- \frac{2}{d}} + 1_{\{d=2\}} N {\log N} \Bigg)^{1/2}.
\end{aligned}
\end{equation}
Note that, since $u \in L^\infty([0,T],C^{1,\alpha}(\mathbb T^d))$, it follows by elliptic regularity that
$$
\| \na u\|_{\L^\infty(\T^d)}+\| \mathfrak{U}\|_{\L^\infty(\T^d)} +\|\na (-\Delta)^{-1} \left( u \,\mathfrak{U}\right)\|_{\C^{0,\alpha}(\T^d)} +\|\pa_t p \|_{\C^{0,\alpha}(\T^d)} +\|\na (-\Delta)^{-1} \left( u \,\mathfrak{U}\right)\|_{\dot\H^{1}(\T^d)} +\|\pa_t p \|_{\dot\H^{1}(\T^d)} \lesssim 1.
$$

Hence
combining \eqref{eq:dt H 2}, \eqref{eq:I1}, \eqref{eq:I2}, and \eqref{eq:I34}, we get 
$$
\frac{\dd}{\dd t} H_{N,\eps} \lesssim H_{N,\eps} + \max \biggl( \frac{1}{\eps^2 N^{\frac{2}{d(d+1)}}},\frac{1}{N^{\frac{\lambda}{d}}}\biggr)+ \eps^2,
$$
hence
\begin{equation}
\label{quanti}
\sup_{t \in [0,T]}H_{N,\eps}(t)\leq e^{CT}\biggl(H_{N,\eps}(0)+\,T\biggl[ \max \biggl( \frac{1}{\eps^2 N^{\frac{2}{d(d+1)}}},\frac{1}{N^{\frac{\lambda}{d}}}\biggr)+ \eps^2\biggr]\biggr).
\end{equation}
Since by assumption $\eps N^{\frac{1}{d(d+1)}} \to +\infty$, the theorem follows.

\section{Proof of Theorem~\ref{thm:2}}
\label{sec3}
The proof of this theorem is very similar to the one of Theorem~\ref{thm:1}, except for the presence of some few extra terms in the time derivative of $\mathscr{H}_{N,\eps}$.
Recalling \eqref{eq:H2},
let us compute the derivative in time of the functional $\mathscr{H}_{N,\eps}$:
\begin{equation}
\label{eq:dt H 2d}
\begin{aligned}
\frac{\dd}{\dd t} \mathscr{H}_{N,\eps}(t)
& =\frac{1}{N} \sum_{i=1}^N  \biggl(u(t,x_i) - \frac{v_i}{\eps}\biggr) \cdot \biggl( \pa_t u(t,x_i) +  \frac{v_i}{\eps} \cdot \na u(t,x_i) + \frac{1}{\eps^2 N} \sum_{\substack{i,j=1\\i\neq j}}^N \na g (x_i -x_j) - \frac{v_i^\perp}{\eps^3}  \biggr) \\
&\quad+    \frac{1}{\eps^3 N^2}\sum_{\substack{i,j=1\\i\neq j}}^N v_i \cdot \na g( x_i - x_j) + \frac{1}{\eps^2} \int_{\Omega \times \Omega \setminus \Delta }  g(x-y)\, \omega(t,x) \,\pa_t \omega(t,y) \, \dd x \,\dd y\\
&\quad- \frac{1}{\eps^3 N}\sum_{i=1}^N \int_{\Omega}  v_i\cdot \na g(x_i -y)  \, \omega(t,y)  \, \dd y - \frac{1}{\eps^2N}\sum_{i=1}^N \int_{\Omega}   g(x_i -y) \, \pa_t \omega(t,y)  \, \dd y\\
&\quad+  \eps^2 \int_{\Omega \times \Omega \setminus \Delta} g(x-y)  \, \mathfrak{U} (t,x) \,\pa_t \mathfrak{U} (t,y) \, \dd x \,\dd y  
\\
&\quad-   \frac{1}{\eps N} \sum_{i=1}^N \int_{\Omega} v_i\cdot \na g (x_i - y)  \, \mathfrak{U} (t,y)\, \dd y -   \frac{1}{N} \sum_{i=1}^N \int_{\Omega}  g (x_i - y) \, \pa_t \mathfrak{U} (t,y)\, \dd y \\
&\quad+   \int_{\Omega \times \Omega \setminus \Delta } g(x-y) \,\pa_t \omega(t,x) \,\mathfrak{U} (t,y) \, \dd x \,\dd y 
+ \int_{\Omega \times \Omega \setminus \Delta } g(x-y)\, \omega(t,x) \,\pa_t \mathfrak{U} (t,y) \, \dd x \,\dd y.
 \end{aligned}
 \end{equation}
Since $u$ solves \eqref{eq:Euler nD}, we have
\begin{equation}
\label{eq:1 2d}
\pa_t u(t,x_i) + \frac{v_i}{\eps} \cdot \na u(t,x_i) = \biggl(\frac{v_i}{\eps}- u(t,x_i)\biggr) \cdot \na u (x_i) - \na p (t,x_i).
\end{equation}
Also, as in the proof of Theorem \ref{thm:1},
we can rewrite (compare with \eqref{eq:2})
$$
\sum_{i=1}^N \int_{\Omega}  v_i\cdot \na g(x_i -y)  \, \omega(t,y)  \, \dd y  = -\sum_{i=1}^N \na \psi (t,x_i) \cdot v_i   =   \sum_{i=1}^N   u^\perp(x_i) \cdot v_i .
$$
Therefore, since $v_i\cdot v_i^\perp=0$ and $u^\perp\cdot v_i=-u\cdot v_i^\perp$, we deduce that
\begin{equation}
\label{eq:2 2d}
-\frac{1}{N} \sum_{i=1}^N  \biggl(u(t,x_i) - \frac{v_i}{\eps}\biggr) \cdot  \frac{v_i^\perp}{\eps^3} 
- \frac{1}{\eps^3 N}\sum_{i=1}^N \int_{\Omega}  v_i\cdot \na g(x_i -y)  \, \omega(t,y)  \, \dd y  =0.
\end{equation}
In addition, we recall that (see  \eqref{eq:2})
\begin{equation}
\label{eq:3 2d}
 \sum_{i=1}^N  \int_{\Omega}  v_i\cdot \na g (x_i - y) \,  \mathfrak{U} (t,y)\, \dd y =  \sum_{i=1}^N v_i  \cdot  \na p (t,x_i ).
\end{equation}
Hence, combining \eqref{eq:1 2d}, \eqref{eq:2 2d}, and \eqref{eq:3 2d}, we can rewrite \eqref{eq:dt H 2d} as
\begin{equation}
\begin{aligned}
\frac{\dd}{\dd t} \mathscr{H}_{N,\eps}(t)
& =
-\frac{1}{N} \sum_{i=1}^N \biggl(u(t,x_i) -\frac{v_i}{\eps}\biggr)\cdot \left[\biggl(u(t,x_i) -\frac{v_i}{\eps}\biggr) \cdot  \na u(t,x_i)\right]\\
&\quad+
\frac{1}{N} \sum_{i=1}^N  u(t,x_i)  \cdot \biggl( -\nabla p(t,x_i) + \frac{1}{\eps^2 N} \sum_{\substack{i,j=1\\i\neq j}}^N \na g (x_i -x_j)   \biggr) \\
&\quad - \frac{1}{\eps^2 N}\sum_{i=1}^N \int_{\Omega}   g(x_i -y) \, \pa_t \omega(t,y)  \, \dd y \\
&\quad\quad\quad\quad + \frac{1}{\eps^2}\int_{\Omega \times \Omega \setminus \Delta }  g(x-y)\, \omega(t,x) \,\pa_t \omega(t,y) \, \dd x \,\dd y+  \int_{\Omega \times \Omega \setminus \Delta } g(x-y) \,\pa_t \omega(t,x) \,\mathfrak{U} (t,y) \, \dd x \,\dd y \\
&\quad 
-\frac{1}{N} \sum_{i=1}^N \int_{\Omega}  g (x_i - y) \, \pa_t \mathfrak{U} (t,y)\, \dd y\\
&\quad\quad \quad\quad +  \int_{\Omega\times \Omega \setminus \Delta}   g (x-y)\, \omega(t,x) \, \pa_t \mathfrak{U} (t,y) \, \dd x \, \dd y +  \eps^2 \int_{\Omega \times \Omega \setminus \Delta} g(x-y) \, \pa_t \mathfrak{U} (t,y) \, \mathfrak{U} (t,x) \, \dd x \,\dd y,
\end{aligned}
\end{equation}
so that
\begin{equation}
\label{eq:dt H 2d bis}
\begin{aligned}
\frac{\dd}{\dd t} \mathscr{H}_{N,\eps}(t)
& =
-\frac{1}{N} \sum_{i=1}^N \biggl(u(t,x_i) -\frac{v_i}{\eps}\biggr)\cdot \left[\biggl(u(t,x_i) -\frac{v_i}{\eps}\biggr) \cdot  \na u(t,x_i)\right]\\
&\quad+
\frac{1}{N} \sum_{i=1}^N  u(t,x_i)  \cdot \biggl( -\nabla p(t,x_i) + \frac{1}{\eps^2 N} \sum_{\substack{i,j=1\\i\neq j}}^N \na g (x_i -x_j)   \biggr) \\
&\quad - \frac{1}{\eps^2} \int_{\Omega}   g(x -y) \, \pa_t \omega(t,y)  \, \dd y\,\dd\biggl(\frac{1}{N}\sum_{i=1}^N \delta_{x_i} -\omega -  \eps^2 \mathfrak{U} \biggr)(x) \\
&\quad- \int_{\Omega}  g (x- y) \, \pa_t \mathfrak{U} (t,y)\, \dd y\,\dd\biggl(\frac{1}{N}\sum_{i=1}^N \delta_{x_i} -\omega -  \eps^2 \mathfrak{U} \biggr)(x). \end{aligned}
 \end{equation}
Now, in analogy with what was done in the proof of Theorem \ref{thm:1},
we observe that
\begin{equation}
\label{eq:4 2d}
\begin{aligned}
&\frac{1}{2\eps^2} \int_{\Omega \times \Omega \setminus \Delta} (u(t,x) - u(t,y)) \cdot \na g (x-y)\, \dd\biggl(  \frac{1}{N}\sum_{j=1}^N \delta_{x_j} -\omega  -   \eps^2 \mathfrak{U} \biggr)(x)\, \dd\biggl(  \frac{1}{N}\sum_{j=1}^N \delta_{x_j} -\omega  -  \eps^2 \mathfrak{U} \biggr)(y) \\
&\quad =   \frac{1}{\eps^2N^2}   \sum_{\substack{i,j=1\\i\neq j}}^N u(t,x_i) \cdot \na g (x_i -x_j) - \frac{1}{N} \sum_{i=1}^N u(t,x_i)  \cdot  \na p (t,x_i )   + \frac{1}{\eps^2N} \sum_{i=1}^N u(t,x_i)  \cdot  \na \psi (t,x_i ) \\
&\quad \quad 
+ \frac{1}{\eps^2} \int_{\Omega}  u(t,y)\cdot \na  g(x -y) \, \omega(t,y)  \, \dd y\,\dd\biggl(\frac{1}{N}\sum_{i=1}^N \delta_{x_i} -\omega -  \eps^2 \mathfrak{U} \biggr)(x)\\
&\quad \quad 
+ \int_{\Omega}  u(t,y)\cdot \na  g(x -y) \, \mathfrak{U} (t,y)  \, \dd y\,\dd\biggl(\frac{1}{N}\sum_{i=1}^N \delta_{x_i} -\omega -  \eps^2 \mathfrak{U} \biggr)(x).
\end{aligned}
\end{equation}
Note first that as $u = \na^\perp \psi$,
\begin{equation}\label{eq:5 2d}
\sum_{i=1}^N u(t,x_i)  \cdot  \na \psi (t,x_i )  = - \sum_{i=1}^N u(t,x_i)  \cdot  u^\perp  (t,x_i ) 
=  0.
\end{equation}
Also, we have
\begin{multline}\label{eq:6 2d}
 \iint_{\Omega\times\Omega\setminus\Delta}  u(t,y)\cdot \na  g(x -y) \, \mathfrak{U} (t,y)  \, \dd y\,\dd\biggl(\frac{1}{N}\sum_{i=1}^N \delta_{x_i} -\omega -  \eps^2 \mathfrak{U} \biggr)(x)\\
  =  \int_{\Omega} \na (-\Delta)^{-1} \left( u \,\mathfrak{U}\right)(x)\, \dd\biggl(\frac{1}{N}\sum_{i=1}^N \delta_{x_i} -\omega -  \eps^2 \mathfrak{U} \biggr)(x).
\end{multline}
Furthermore, using \eqref{eq:Euler 2D} and  \eqref{eq:p U}, 
we see that
\begin{multline}\label{eq:7 2d}
 \iint_{\Omega\times\Omega\setminus\Delta}u(t,y)\cdot \na  g(x -y) \, \omega(t,y)  \, \dd y\,\dd\biggl(\frac{1}{N}\sum_{i=1}^N \delta_{x_i} -\omega -  \eps^2 \mathfrak{U} \biggr)(x)\\
=- \int_{\Omega}  g(x -y) \,\pa_t \omega(t,y)  \, \dd y\,\dd\biggl(\frac{1}{N}\sum_{i=1}^N \delta_{x_i} -\omega -  \eps^2 \mathfrak{U} \biggr)(x),
\end{multline}
and
\begin{multline}\label{eq:8 2d}
  \iint_{\Omega\times\Omega\setminus\Delta}g(x -y) \, \partial_t\mathfrak{U} (t,y)  \, \dd y\,\dd\biggl(\frac{1}{N}\sum_{i=1}^N \delta_{x_i} -\omega -  \eps^2 \mathfrak{U} \biggr)(x)\\
= \int_{\Omega}  \pa_t p (t,x)\, \dd\biggl(\frac{1}{N}\sum_{i=1}^N \delta_{x_i} -\omega -  \eps^2 \mathfrak{U} \biggr)(x).
\end{multline}
Combining together \eqref{eq:dt H 2d bis}, \eqref{eq:4 2d}, \eqref{eq:5 2d}, \eqref{eq:6 2d}, \eqref{eq:7 2d}, and  \eqref{eq:8 2d}, we end up with
\begin{equation}
\begin{aligned}
\frac{\dd}{\dd t} \mathscr{H}_{N,\eps} &= -\frac{1}{N} \sum_{i=1}^N \biggl(u(t,x_i) -\frac{v_i}{\eps}\biggr)\cdot \left[\biggl(u(t,x_i) -\frac{v_i}{\eps}\biggr) \cdot  \na u(t,x_i)\right]   \\
&+ \frac{1}{2\eps^2} \int_{\Omega \times \Omega \setminus \Delta} (u(t,x) - u(t,y)) \cdot \na g (x-y)\, \dd\biggl(  \frac{1}{N}\sum_{i=1}^N \delta_{x_i} - \omega- \eps^2 \mathfrak{U} \biggr)(x)\, \dd\biggl(  \frac{1}{N}\sum_{j=1}^N \delta_{x_j} - \omega-   \eps^2 \mathfrak{U} \biggr)(y) \\
&-  \int_{\Omega} \na (-\Delta)^{-1} \left( u \,\mathfrak{U}\right)(x)\, \dd\biggl(\frac{1}{N}\sum_{i=1}^N \delta_{x_i} -\omega -  \eps^2 \mathfrak{U} \biggr)(x)  \\
&+\int_{\Omega}  \pa_t p (t,x)\, \dd\biggl(\frac{1}{N}\sum_{i=1}^N \delta_{x_i} -\omega -  \eps^2 \mathfrak{U} \biggr)(x) =: \sum_{i=1}^4 J_i.
\end{aligned}
\end{equation}
Each term $J_i$ can now be bounded exactly as was done in the proof of Theorem~\ref{thm:1} to estimate the terms $I_i$,
replacing $1+ \eps^2 \| \mathfrak{U}\|_{\L^\infty(\Omega)} $ by $ \| \omega\|_{\L^\infty(\Omega)} + \eps^2 \| \mathfrak{U}\|_{\L^\infty(\Omega)} $ so that
$$
\frac{\dd}{\dd t} \mathscr H_{N,\eps} \lesssim \mathscr H_{N,\eps} + \max \biggl( \frac{1}{\eps^2 N^{\frac{1}{3}}},\frac{1}{N^{\frac{\lambda}{2}}}\biggr)+ \eps^2.
$$
Therefore
$$
\sup_{t \in [0,T]}\mathscr H_{N,\eps}(t)\leq e^{CT}\biggl(\mathscr H_{N,\eps}(0)+\,T\biggl[ \max \biggl( \frac{1}{\eps^2 N^{\frac{1}{3}}},\frac{1}{N^{\frac{\lambda}{2}}}\biggr)+ \eps^2\biggr]\biggr).
$$
Since by assumption $\eps N^{\frac{1}{6}} \to +\infty$, the theorem follows.

\bibliography{Mikabib.bib}
\bibliographystyle{abbrv}

%
\end{document}